\documentclass[12pt]{article}
\usepackage{amsmath,amssymb,amsthm,booktabs}
\usepackage{mathrsfs}
\usepackage{geometry}
\usepackage{pifont}
\usepackage[OT1]{fontenc}
\usepackage[utf8]{inputenc}
\usepackage[colorlinks,citecolor=blue,urlcolor=blue]{hyperref}
\usepackage{txfonts}
\usepackage{indentfirst}
\usepackage{multirow}
\usepackage{float,subfig}

\usepackage{bm}
\usepackage{euscript}
\usepackage{graphicx}
\usepackage{multicol}
\usepackage[usenames,dvipsnames,svgnames,table]{xcolor}

\usepackage[numbers]{natbib}
\numberwithin{equation}{section} \theoremstyle{plain}
\newtheorem{theorem}{Theorem}[section]
\newtheorem{lemma}{Lemma}[section]
\newtheorem{corollary}{Corollary}[section]

\newtheorem{example}{Example}[section]

\newtheorem{remark}{Remark}[section]

\begin{document}

\newcommand{\gai}[1]{{#1}}


\makeatletter
\def\ps@pprintTitle{%
  \let\@oddhead\@empty
  \let\@evenhead\@empty
  \let\@oddfoot\@empty
  \let\@evenfoot\@oddfoot
}
\makeatother

\newcommand\tabfig[1]{\vskip5mm \centerline{\textsc{Insert #1 around here}}  \vskip5mm}

\vskip2cm

\title{$L^p$ estimations of fully coupled FBSDEs}
\author{Qingxin Meng\thanks{School of Mathematical Sciences, Huzhou University, Huzhou 313000,
China, Email: mqx@zjhu.edu.cn.This work was supported by  National Natural Science Foundation of China (No. 11871121), and by the Natural Science Foundation of Zhejiang Province (No. LY21A010001).} \quad  Shuzhen Yang\thanks{Shandong University-Zhong Tai Securities Institute for Financial Studies, Shandong University, PR China, (yangsz@sdu.edu.cn). This work was supported by the National Key R\&D program of China (Grant No.2018YFA0703900), National Natural Science Foundation of China (Grant No.11701330), and Young Scholars Program of Shandong University.}
}
\date{}
\maketitle

\begin{abstract}
 In this study, for any given terminal time $T$, we establish an $L^p$ ($P>2$) estimations of fully coupled FBSDEs based on the $L^2$ estimations. \citet{Y20} proposed that a natural question is whether an adapted $L^2$-solution is an adapted $L^p$-solution for some $p>2$. In this study, we give a positive answer to this question. For any given terminal time $T$, based on an observation of the relation between $L^2$ and $L^p$ estimations of FBSDEs, we prove that a unique $L^2$-solution of fully coupled FBSDEs is an $L^p$-solution under standard conditions on the coefficients. Furthermore, we show that the fully coupled FBSDEs developed in the linear quadratic optimal control problem or investigated by the decoupling random field method admit a unique $L^p$-solution.
\end{abstract}



\section{Introduction}

Since \citet{Bismut} introduced the linear backward stochastic differential equations (BSDEs) from the stochastic optimal control problem, and \citet{PardouxE3} established the existence and uniqueness theorem for nonlinear BSDEs. Many related studies on BSDEs have been conducted.
\citet{PengS5} and  \citet{PardouxE2} generalized the classical Feynman-Kac formula by establishing the relationship between BSDEs and systems of quasilinear parabolic PDEs. More details see \citet{PengS1}. The application of BSDEs in finance was concluded in \citet{KPC97}.

Subsequently, many authors have begun studying the fully coupled forward-backward SDEs (FBSDEs), which are related to the  optimal control problem and a new quasilinear parabolic PDE. The fixed-point method, method of continuation, and decoupling random field are the three main methods used to study the existence and uniqueness of an adapted $L^2$-solution of fully coupled FBSDEs. \citet{An93} investigated the existence and uniqueness results of fully coupled FBSDEs within a sufficiently small terminal time $T$. \citet{PT99} extended these results based on a fixed-point approach. Using a decoupling random field method, \citet{MaJ} originally observed the well-known Four Step Scheme for fully coupled FBSDEs. Further results see \citet{MY99,De02,Zh97,MYZ,MW15}. The monotonicity conditions and also called the method of continuation developed in  \citet{HP95}, and see \citet{PW99,Yo97}. Additionally, the monographs \citet{MY99} and \citet{CZ13} for the theories of FBSDEs and applications in finance.

The $L^p$ estimations of fully coupled FBSDEs play an important role in theories and related applications. For example, stochastic optimal control problems for fully coupled FBSDEs can be developed. More details are available in \citet{Y10,LW14,LW141,SLY16,Hu17,Hu18} for further details. However, there are many difficulties in deriving the $L^p$ estimations for fully coupled FBSDEs. Currently, the $L^p$ estimations of fully coupled FBSDEs are established based on some strong assumptions on the coefficients of FBSDEs and within a sufficiently small $T$ (\citet{MY99} and \citet{CZ13}). \citet{Y19} established an $L^p$-Theory for fully FBSDEs with random coefficients on small durations when the Lipschatiz constant of $Z$ in diffusion term does not to be small.

In the field of fully coupled FBSDEs, $L^p$ ($p>2$) estimations of fully coupled FBSDEs are  different from $L^2$ estimations, where the $L^2$ estimations are derived by the continuation method and decoupling random field method for any given terminal time $T$. One can only obtain a $L^p$ ($p>2$) estimations for a sufficiently small $T$ using the fixed-point method. Thus, for any given terminal time $T$, obtaining $L^p$ ($p>2$) estimations of general fully coupled FBSDEs remains an open problem. In one-dimensional case, \citet{MW15} developed a uniform approach for fully coupled FBSDEs by decoupling random field method, and established an $L^p$ ($p>2$) estimations when the diffusion term does not depend on solution $Z$. Recently, \citet{Y20} proposed the question of whether an adapted $L^2$-solution is an adapted $L^P$-solution for some $p>2$. Following the question in \citet{Y20}, we provide a positive answer to the open problem. For a given terminal time $T>0$, we prove that a unique $L^2$-solution of fully coupled FBSDEs is an $L^p$-solution under the usual Lipschatiz ( Lipschatiz constant of $Z$ in diffusion term should be sufficiently small) and linear growth conditions on the coefficients. Then, we apply our results to study fully coupled liner FBSDEs which are generalized by a linear quadratic optimal control problem, and obtain the $L^p$ estimations for fully coupled liner FBSDEs with random coefficients. Based on the results of the decoupling random field method in fully coupled FBSDEs, under the same conditions, we improve the $L^2$ estimations to the $L^p$ estimations for any given terminal time $T$.

The contributions of this study are given as follows:

i). Based on an $L^2$ estimations of fully couple FBSDEs, we derive a uniform estimations for the solution of the backward SDE on the solution of the forward SDE;

ii). Applying the uniform estimations, we establish an $L^p$ estimations for fully couple FBSDEs;

iii). Under the standard Lipschatiz ( Lipschatiz constant of $Z$ in diffusion term should be sufficiently small) and linear growth conditions, we claim that the $L^2$ estimations of fully couple FBSDEs are equivalent to $L^p$ estimations;

iv). In one-dimensional case, \citet{MW15} established an $L^p$ ($p>2$) estimations by decoupling random field method. Different from the method in \citet{MW15}, we use an $L^2$ estimations to establish the $L^p$ estimations for multi-dimensional case. Thus, our results can improve the wellposedness of fully coupled FBSDEs in $L^2$ to $L^p$.

The remainder of this study is organized as follows. In Section \ref{sec:2}, we introduce the notations, basic assumptions, and preliminary results. Then, we establish the main results that is the $L^2$-solution of a fully coupled FBSDE is an $L^p$-solution. In Section \ref{sec:3}, we consider a linear-quadratic optimal control problem and prove that the fully coupled linear FBSDE admits the $L^p$-solution, especially for the random coefficient case. Based on the "decoupling random field"  method, we show that the $L^2$-solution of a fully coupled FBSDE is an $L^p$-solution in Section \ref{sec:4}. We then extend the main results of this study in Section \ref{sec:4.5}. Finally, we conclude the study in Section \ref{sec:5}.

\section{Main results}\label{sec:2}

Let $B$ be a one-dimensional standard Brownian motion defined on a complete
filtered probability space $(\Omega,\mathcal{F},P;\{ \mathcal{F}(t)\}_{t\geq
0})$, where $\{ \mathcal{F}(t)\}_{t\geq0}$ is the $P$-augmentation of the
natural filtration generated by  $B$. We consider the following fully coupled FBSDE, parameterized by the initial condition $(t,\xi) \in[0,T] \times L^{2}(\Omega,\mathcal{F}%
_{t},P;\mathbb{R}^{n})$,
\begin{equation}
\label{equ2.1}\left \{
\begin{array}
[c]{llll}%
\mathrm{d}X_{s}^{t,\xi} & = & b(s,\Theta_{s}^{t,\xi})\mathrm{d}s + \sigma
(s,\Theta_{s}^{t,\xi}) \mathrm{d}B_{s}, & \\
\mathrm{d}Y_{s}^{t,\xi} & = & -f(s,\Theta_{s}^{t,\xi})\mathrm{d}s +
Z_{s}^{t,\xi}\mathrm{d}B_{s}, \  \  \  \  \ s\in[t,T], & \\
X_{t}^{t,\xi} & = & \xi, & \\
Y_{T}^{t,\xi} & = & \Phi(X_{T}^{t,\xi}), &
\end{array}
\right.
\end{equation}
where $\Theta_s^{t,\xi}=(X_s^{t,\xi},Y_s^{t,\xi},Z_s^{t,\xi}),\ t\leq s\leq T$, and
$$
\begin{array}
[c]{ll}
& \! \! \! \! \! b: \Omega \times[0,T] \times \mathbb{R}^{n} \times \mathbb{R}^m\times \mathbb{R}^{m}
 \longrightarrow \mathbb{R}^{n},\\
& \! \! \! \! \!  \sigma: \Omega \times[0,T] \times \mathbb{R}^{n}
\times \mathbb{R}^m\times \mathbb{R}^{m} \longrightarrow
\mathbb{R}^{n},\\
& \! \! \! \! \!  f: \Omega \times[0,T] \times \mathbb{R}^{n} \times \mathbb{R}^m
\times \mathbb{R}^{m} \longrightarrow \mathbb{R}^m,\\
&\! \! \! \! \! \Phi: \Omega \times \mathbb{R}^n \longrightarrow\mathbb{R}^m,
\end{array}
$$
satisfy the following linear growth, Lipschatiz continuous assumptions.
\begin{description}
\item[$( \mathbf{H2.1})$] For any $t\in[0,T]$, and $(x,y,z)\in
\mathbb{R}^{n}\times \mathbb{R}^m\times \mathbb{R}^{m},\  \mbox{P-a.s.},$ there exists $L$ such that
$$
 |b(t,x,y,z)|+|\sigma(t,x,y,z)|+ |f(t,x,y,z)|+|\Phi(x)| \leq L(1+|x|+|y|+|z|). $$
\end{description}
\begin{description}
\item[$( \mathbf{H2.2})$] There exist constant $K$ and
a sufficiently small constant $L_{\sigma}\geq0$ such that for all $t\in[0,T],\ x_{1},x_{2}%
\in \mathbb{R}^{n},\ y_{1},y_{2}\in \mathbb{R}^m,\ z_{1},z_{2}\in \mathbb{R}%
^{m}$,
\[
|b(t,x_{1},y_{1},z_{1})-b(t,x_{2},y_{2},z_{2})|\leq
K(|x_{1}-x_{2}|+|y_{1}-y_{2}|+|z_{1}-z_{2}|),
\]
\[
|\sigma(t,x_{1},y_{1},z_{1})-\sigma(t,x_{2},y_{2},z_{2})|\leq
K(|x_{1}-x_{2}|+|y_{1}-y_{2}|)+L_{\sigma}|z_{1}-z_{2}|,
\]
\[
|f(t,x_{1},y_{1},z_{1})-f(t,x_{2},y_{2},z_{2})|\leq
K(|x_{1}-x_{2}|+|y_{1}-y_{2}|+|z_{1}-z_{2}|),
\]
\[
|\Phi(x_{1})-\Phi(x_{2})|\leq
K|x_{1}-x_{2}|.
\]
\end{description}

Based on Assumptions $( \mathbf{H2.1})$ and $( \mathbf{H2.2})$, we obtain the following classical results, which can be found in \citet{An93,De02,LW14,LW141,Y20}.
\begin{lemma}
\label{lem-1} Let Assumptions $( \mathbf{H2.1}),\ ( \mathbf{H2.2})$ hold. Then, there exists a constant $\delta>0$ depending on $(K,L_\sigma)$, such that FBSDE
(\ref{equ2.1}) admits a unique solution
$\Theta_{s}^{t,\xi}$ on a
interval $[t, t+\delta]$, where $t+\delta=T$.
\end{lemma}

\begin{lemma}
\label{lem-2} Let Assumptions $( \mathbf{H2.1}),\ (\mathbf{H2.2})$
hold. Then, for any given $p\geq 2,$ there exists a constant ${\delta} > 0$ that depends on $(K,L_\sigma)$, and constant ${C}_0$ that depends on $(p,L,K,L_\sigma)$ such that
for every $\xi \in L^{p}(\Omega,\mathcal{F}_{t}, P;\mathbb{R}^{n}),$
\[%
\begin{array}
[c]{llll}%
 &  E [ \mathop{\rm sup}\limits_{t\leq s\leq t+\delta}
|X^{t,\xi}_{s} |^{p} + \mathop{\rm sup}\limits_{t\leq s\leq
t+\delta} |Y^{t,\xi}_{s} |^{p} + (\int_{t}^{t+\delta}
|Z^{t,\xi}_{s}|^{2}\mathrm{d}s)^{\frac{p}{2}} \mid \mathcal{F}_{t}]\leq {C}_0(1 + |\xi|^{p}).& \\
\end{array}
\]
\end{lemma}

\subsection{$L^p$ estimations}

We first consider an example which motivates this study.
\begin{example}
Let $P_s=\int_t^sa_r\mathrm{d}r$, and let $a(\cdot),b(\cdot),c(\cdot)$ be a real-valued one-dimensional bounded progressively measurable process on $[t,T]$.
\begin{equation}
\label{ex-equ1}\left \{
\begin{array}
[c]{llll}%
\mathrm{d}X^{t,\xi}_{s} & = & b_sY^{t,\xi}_s\mathrm{d}s +c_s \mathrm{d}B_{s}, & \\
\mathrm{d}Y^{t,\xi}_{s}& = & [a_sX^{t,\xi}_s+b_sP_sY^{t,\xi}_s]\mathrm{d}s +
Z^{t,\xi}_{s}\mathrm{d}B_{s}, \  \  \  \  \ s\in[t,T], & \\
X^{t,\xi}_{t}& = & \xi, & \\
Y^{t,\xi}_{T} & = & P_TX^{t,\xi}_T. &
\end{array}
\right.
\end{equation}
It is easy to verify that $Y^{t,\xi}_{s}=P_sX^{t,\xi}_s$ satisfies the backward equation of (\ref{ex-equ1}). Then, the forward equation can be rewritten as follows:
\begin{equation}\label{ex-equ2}
\mathrm{d}X^{t,\xi}_{s} = b_sP_sX^{t,\xi}_s\mathrm{d}s +c_s \mathrm{d}B_{s},\ X^{t,\xi}_{t}=\xi.
\end{equation}

In the following, we assume that Equation (\ref{ex-equ1}) admits a unique $L^2$-solution. From $Y^{t,\xi}_{s}=P_sX^{t,\xi}_s$, we have the following:
\begin{equation}\label{ex-equ3}
\left|Y^{t,\xi}_{s}-Y^{t,\xi'}_{s}\right|\leq C\left|X^{t,\xi}_{s}-X^{t,\xi'}_{s}\right|,\ t\leq s\leq T.
\end{equation}

From Equation (\ref{ex-equ2}), we can obtain the $L^p$ estimations for $\sup_{t\leq s\leq T}\left|X^{t,\xi}_{s}-X^{t,\xi'}_{s}\right|$, that is,
$$
E[\sup_{t\leq s\leq T}\left|X^{t,\xi}_{s}-X^{t,\xi'}_{s}\right|^p]\leq C_1\left|\xi-\xi'\right|^p.
$$
Combining with Equation (\ref{ex-equ3}), we obtain the $L^p$ ($p>2$) estimations for $\sup_{t\leq s\leq T}\left|Y^{t,\xi}_{s}-Y^{t,\xi'}_{s}\right|$,
$$
E[\sup_{t\leq s\leq T}\left|Y^{t,\xi}_{s}-Y^{t,\xi'}_{s}\right|^p]\leq C_1\left|\xi-\xi'\right|^p.
$$

Then, using the backward equation of (\ref{ex-equ1}), we have the following $L^p$ estimations:
$$
E [ \mathop{\rm sup}\limits_{t\leq s\leq T}
|X^{t,\xi}_{s}-X^{t,\xi'}_{s} |^{p} + \mathop{\rm sup}\limits_{t\leq s\leq T} |Y^{t,\xi}_{s}-Y^{t,\xi'}_{s} |^{p} + \big(\int_{t}^{T}
|Z^{t,\xi}_{s}-Z^{t,\xi'}_{s}|^{2}\mathrm{d}s\big)^{\frac{p}{2}}]\leq {C}_1|\xi-\xi'|^{p}.
$$
We can see that  (\ref{ex-equ2}) is a useful condition for proving the general $L^p$ estimations of some fully coupled FBSDE. Thus, in the following, we use an $L^2$ estimations of fully coupled FBSDE to obtain the condition (\ref{ex-equ3}).
\end{example}

In the following, we establish the equivalence between the $L^2$ and $L^p$ ($p>2$) estimations of FBSDE (\ref{equ2.1}) under Assumptions $(\mathbf{H2.1}),\ (\mathbf{H2.2})$.
\begin{theorem}
\label{the-1} Let Assumptions $( \mathbf{H2.1}),\ (\mathbf{H2.2})$
hold, and we assume that for every $\xi,\xi' \in L^{p}(\Omega,\mathcal{F}_{t}, P;\mathbb{R}^{n})$, $L^2$ estimations of FBSDE (\ref{equ2.1}) are right,
\begin{equation}
\label{equa2.2}
 E [ \mathop{\rm sup}\limits_{t\leq s\leq T}
|X^{t,\xi}_{s} |^{2} + \mathop{\rm sup}\limits_{t\leq s\leq T} |Y^{t,\xi}_{s} |^{2} + (\int_{t}^{T}
|Z^{t,\xi}_{s}|^{2}\mathrm{d}s)\mid \mathcal{F}_{t} ]\leq {C}_1(1 + |\xi|^{2}),
\end{equation}
\begin{equation}
\label{equa2.3}
E [ \mathop{\rm sup}\limits_{t\leq s\leq T}
|X^{t,\xi}_{s}-X^{t,\xi'}_{s} |^{2} + \mathop{\rm sup}\limits_{t\leq s\leq T} |Y^{t,\xi}_{s}-Y^{t,\xi'}_{s} |^{2} + \int_{t}^{T}
|Z^{t,\xi}_{s}-Z^{t,\xi'}_{s}|^{2}\mathrm{d}s\mid \mathcal{F}_{t} ]\leq {C}_1|\xi-\xi'|^{2},
\end{equation}
where $C_1$ is a positive constant and independent from $t\in [0,T]$.

Then, FBSDE (\ref{equ2.1}) admits a unique $L^p$ ($p>2$) solution with $t=0$ and any given terminal time $T$, that is,
\[%
\begin{array}
[c]{llll}%
 &  E [ \mathop{\rm sup}\limits_{0\leq s\leq T}
|X^{0,\xi}_{s} |^{p} + \mathop{\rm sup}\limits_{0\leq s\leq T} |Y^{0,\xi}_{s} |^{p} + (\int_{0}^{T}
|Z^{0,\xi}_{s}|^{2}\mathrm{d}s)^{\frac{p}{2}}]\leq {C}_2(1 + |\xi|^{p}),& \\
 &  E [ \mathop{\rm sup}\limits_{0\leq s\leq T}
|X^{0,\xi}_{s}-X^{0,\xi'}_{s} |^{p} + \mathop{\rm sup}\limits_{0\leq s\leq T} |Y^{0,\xi}_{s}-Y^{0,\xi'}_{s} |^{p} + (\int_{0}^{T}
|Z^{0,\xi}_{s}-Z^{0,\xi'}_{s}|^{2}\mathrm{d}s)^{\frac{p}{2}}]\leq {C}_2|\xi-\xi'|^{p}.& \\
\end{array}
\]
\end{theorem}

\noindent \textbf{Proof}:  Based on $L^2$ estimations (\ref{equa2.2}) and (\ref{equa2.3}), FBSDE (\ref{equ2.1}) admits a unique $L^2$-solution in the interval $[t,T]$ for any given $(t,\xi)\in [0,T]\times L^2(\Omega,\mathcal{F}_t,P; \mathbb{R}^n)$.

We consider FBSDE (\ref{equ2.1}) at the initial time $0$:
\begin{equation}
\label{equ2.4}\left \{
\begin{array}
[c]{llll}%
dX_{s}^{0,\xi} & = & b(s,\Theta_{s}^{0,\xi})\mathrm{d}s + \sigma
(s,\Theta_{s}^{0,\xi}) \mathrm{d}B_{s}, & \\
dY_{s}^{0,\xi} & = & -f(s,\Theta_{s}^{0,\xi})\mathrm{d}s +
Z_{s}^{0,\xi}\mathrm{d}B_{s}, \  \  \  \  \ s\in[0,T], & \\
Y_{T}^{0,\xi} & = & \Phi(X_{T}^{0,\xi}). &
\end{array}
\right.
\end{equation}
Using inequality (\ref{equa2.3}), we have the following:
$$
E [ \mathop{\rm sup}\limits_{t\leq s\leq T}
|X^{0,\xi}_{s}-X^{0,\xi'}_{s} |^{2} + \mathop{\rm sup}\limits_{t\leq s\leq T} |Y^{0,\xi}_{s}-Y^{0,\xi'}_{s} |^{2} + \int_{t}^{T}
|Z^{0,\xi}_{s}-Z^{0,\xi'}_{s}|^{2}\mathrm{d}s\mid \mathcal{F}_{t} ]\leq {C}_1|X^{0,\xi}_{t}-X^{0,\xi'}_{t} |^{2},
$$
and thus
$$
\left|Y_t^{0,\xi}-Y_t^{0,\xi'}\right|^2\leq E[\sup_{t\leq r\leq T}\left|Y_r^{0,\xi}-Y_r^{0,\xi'}\right|^2 \mid \mathcal{F}_{t} ] \leq C_1 \left| X_t^{0,\xi}- X_t^{0,\xi'}\right|^2,\quad t\in [0,T].
$$
Thus, the solution $(Y_t^{0,\xi},Y_t^{0,\xi'})$ satisfies
\begin{equation}\label{eq001}
\left|Y_t^{0,\xi}-Y_t^{0,\xi'}\right|\leq \sqrt{C_1} \left| X_t^{0,\xi}- X_t^{0,\xi'}\right|,\quad t\in [0,T].
\end{equation}

Now, we consider FBSDE (\ref{equ2.1}) in the interval $[0,\delta]$, where $\delta>0$ is a constant which be given later.
Based on Assumptions $( \mathbf{H2.1})$, and $(\mathbf{H2.2})$ and inequality (\ref{eq001}), using Lemma \ref{lem-2}, there exists a constant $\delta>0$ such that
\[%
\begin{array}
[c]{llll}%
 & & E[ \mathop{\rm sup}\limits_{0\leq s\leq \delta}
|X^{0,\xi}_{s} |^{p}+ \mathop{\rm sup}\limits_{0\leq s\leq \delta}|Y^{0,\xi}_{s} |^{p}+ (\int_{0}^{\delta}
|Z^{0,\xi}_{s}|^{2}ds)^{\frac{p}{2}} \mid \mathcal{F}_{0}]\leq {C}^{(1)}_1(1 + |X^{0,\xi}_{0}|^{p}), \\
\end{array}
\]
where $C^{(1)}_1$ depends on the constants ${C_1}$ and $L,K,L_{\sigma}$ in Assumptions $( \mathbf{H2.1}),\ (\mathbf{H2.2})$. Note that, the coefficients of FBSDE (\ref{equ2.1}) satisfy the same assumptions in the interval $[0,T]$. Then, combining inequality (\ref{eq001}), we can obtain the following inequality by inductive method for the same $\delta$,
\[%
\begin{array}
[c]{llll}%
 & & E[ \mathop{\rm sup}\limits_{(i-1)\delta\leq s\leq i\delta}
|X^{0,\xi}_{s} |^{p}+ \mathop{\rm sup}\limits_{(i-1)\delta\leq s\leq i\delta}|Y^{0,\xi}_{s} |^{p}+ (\int_{(i-1)\delta}^{i\delta}
|Z^{0,\xi}_{s}|^{2}ds)^{\frac{p}{2}} \mid \mathcal{F}_{(i-1)\delta}]\leq {C}^{(1)}_1(1 + |X^{0,\xi}_{(i-1)\delta}|^{p}), \\
\end{array}
\]
where $1\leq i\leq k$, and $k$ is a positive integer. Without loss of generality, we assume that $T=k\delta$.

We first consider the cases where $i=1,2$,
\[%
\begin{array}
[c]{llll}%
 & & E[ \mathop{\rm sup}\limits_{0\leq s\leq \delta}
|X^{0,\xi}_{s} |^{p}+ \mathop{\rm sup}\limits_{0\leq s\leq \delta}|Y^{0,\xi}_{s} |^{p}+ (\int_{0}^{\delta}
|Z^{0,\xi}_{s}|^{2}ds)^{\frac{p}{2}} \mid \mathcal{F}_{0}]\leq C^{(1)}_1(1 + |\xi|^{p}), & \\
\end{array}
\]
and
\[%
\begin{array}
[c]{llll}%
 & & E[ \mathop{\rm sup}\limits_{\delta\leq s\leq 2\delta}
|X^{0,\xi}_{s} |^{p}+ \mathop{\rm sup}\limits_{\delta\leq s\leq 2\delta} |Y^{0,\xi}_{s} |^{p}+ (\int_{\delta}^{2\delta}
|Z^{0,\xi}_{s}|^{2}ds)^{\frac{p}{2}} \mid \mathcal{F}_{\delta}]\leq C^{(1)}_1(1 + |X^{0,\xi}_{\delta}|^{p}). & \\
\end{array}
\]
From the case $i=1$, we have
$$
E[{C}^{(1)}_{1}(1 + |X^{0,\xi}_{\delta}|^{p}) \mid \mathcal{F}_{0}]\leq C^{(1)}_1\left(1+C^{(1)}_1(1 + |\xi|^{p})\right)\leq (C^{(1)}_1+(C^{(1)}_1)^2)(1 + |\xi|^{p}).
$$
Let ${C}^{(2)}_{1}=2{C}^{(1)}_{1}+({C}^{(1)}_{1})^2$, it follows that
$$
{C}^{(1)}_{1}(1 + |\xi|^{p})+E[{C}^{(1)}_{1}(1 + |X^{0,\xi}_{\delta}|^{p}) \mid \mathcal{F}_{0}]\leq {C}^{(2)}_{1}(1 + |\xi|^{p}).
$$
Adding on both sides of cases $i=1$ and $i=2$, we have
\[%
\begin{array}
[c]{llll}%
 & & E[ \mathop{\rm sup}\limits_{0\leq s\leq 2\delta}
|X^{0,\xi}_{s} |^{p}+ \mathop{\rm sup}\limits_{0\leq s\leq 2\delta}|Y^{0,\xi}_{s} |^{p}+(\int_{0}^{\delta}
|Z^{0,\xi}_{s}|^{2}ds)^{\frac{p}{2}}+ (\int_{\delta}^{2\delta}
|Z^{0,\xi}_{s}|^{2}ds)^{\frac{p}{2}} \mid \mathcal{F}_{0}]\leq {C}^{(2)}_{1}(1 + |\xi|^{p}). & \\
\end{array}
\]
Combining the inequality,
$$
(a+b)^k\leq 2^k(a^k+b^k),\ 0\leq a,b,\ 1\leq k.
$$
Let $\hat{C}^{(2)}_{1}=2^{\frac{p}{2}}{C}^{(2)}_{1}$, one obtains
\[%
\begin{array}
[c]{llll}%
 & & E[ \mathop{\rm sup}\limits_{0\leq s\leq 2\delta}
|X^{0,\xi}_{s} |^{p}+ \mathop{\rm sup}\limits_{0\leq s\leq 2\delta}|Y^{0,\xi}_{s} |^{p}+ (\int_{0}^{2\delta}
|Z^{0,\xi}_{s}|^{2}ds)^{\frac{p}{2}} \mid \mathcal{F}_{0}]\leq \hat{C}^{(2)}_{1}(1 + |\xi|^{p}). & \\
\end{array}
\]

Then, we consider the case $i=3$,
\[%
\begin{array}
[c]{llll}%
 & & E[ \mathop{\rm sup}\limits_{2\delta\leq s\leq 3\delta}
|X^{0,\xi}_{s} |^{p}+ \mathop{\rm sup}\limits_{2\delta\leq s\leq 3\delta}|Y^{0,\xi}_{s} |^{p}+ (\int_{2\delta}^{3\delta}
|Z^{0,\xi}_{s}|^{2}ds)^{\frac{p}{2}} \mid \mathcal{F}_{2\delta}]\leq {C}^{(1)}%
_{1}(1 + |X^{0,\xi}_{2\delta}|^{p}). & \\
\end{array}
\]
Similar with the above analysis, we have
\[%
\begin{array}
[c]{llll}%
 & & E[ \mathop{\rm sup}\limits_{0\leq s\leq 3\delta}
|X^{0,\xi}_{s} |^{p}+ \mathop{\rm sup}\limits_{0\leq s\leq 3\delta}|Y^{0,\xi}_{s} |^{p}+ (\int_{0}^{3\delta}
|Z^{0,\xi}_{s}|^{2}ds)^{\frac{p}{2}} \mid \mathcal{F}_{0}]\leq \hat{C}^{(3)}_{1}(1 + |\xi|^{p}). & \\
\end{array}
\]
Then, considering the case $i=4,\cdots,k$, we can obtain the $L^p$ estimations
\[%
\begin{array}
[c]{llll}%
 & & E[ \mathop{\rm sup}\limits_{0\leq s\leq k\delta}
|X^{0,\xi}_{s} |^{p}+\mathop{\rm sup}\limits_{0\leq s\leq k\delta} |Y^{0,\xi}_{s} |^{p}+ (\int_{0}^{k\delta}
|Z^{0,\xi}_{s}|^{2}ds)^{\frac{p}{2}} \mid \mathcal{F}_{0}]\leq \hat{C}^{(k)}_{1}(1 + |\xi|^{p}),& \\
\end{array}
\]
where $k\delta=T$.

Now, let $C_2=\hat{C}^{(k)}_{1}$. Then, FBSDE (\ref{equ2.1}) admits a unique $L^p$-solution with $t=0$. This completes the proof. $\quad  \qquad  \Box$

\begin{remark}
Theorem \ref{the-1} establishes the $L^p$ estimations for fully coupled FBSDE with any given terminal time $T$. From the proof of Theorem \ref{the-1}, we can see that the basic Lipschatiz and linear growth conditions on $b,\sigma,f,\Phi$ and  sufficiently small Lipschatiz constant $L_{\sigma}$ of the diffusion term $\sigma$ on $Z$ are necessary. Furthermore, we obtain the following inequality
$$
\left|Y_s^{0,\xi}-Y_s^{0,\xi'}\right|\leq \sqrt{C_1} \left| X_s^{0,\xi}- X_s^{0,\xi'}\right|,\quad s\in [0,T],
$$
from the $L^2$ conditional expectation estimations. Thus, we can establish the following corollary from Theorem \ref{the-1}.
\end{remark}

\begin{corollary}

\label{coro-1-1} Let Assumptions $( \mathbf{H2.1}),\ (\mathbf{H2.2})$
hold, and we assume that for every $\xi,\xi' \in L^{p}(\Omega,\mathcal{F}_{t}, P;\mathbb{R}^{n})$, there exists a positive constant $C_1>0$ such that
$$
\left|Y_s^{0,\xi}-Y_s^{0,\xi'}\right|\leq {C_1}\left| X_s^{0,\xi}- X_s^{0,\xi'}\right|,\quad s\in [0,T].
$$

Then, FBSDE (\ref{equ2.1}) admits a unique $L^p$ ($p>2$) solution with $t=0$, and
\[%
\begin{array}
[c]{llll}%
 &  E [ \mathop{\rm sup}\limits_{0\leq s\leq T}
|X^{0,\xi}_{s} |^{p} + \mathop{\rm sup}\limits_{0\leq s\leq T} |Y^{0,\xi}_{s} |^{p} + (\int_{0}^{T}
|Z^{0,\xi}_{s}|^{2}\mathrm{d}s)^{\frac{p}{2}}]\leq {C}_2(1 + |\xi|^{p}),& \\
 &  E [ \mathop{\rm sup}\limits_{0\leq s\leq T}
|X^{0,\xi}_{s}-X^{0,\xi'}_{s} |^{p} + \mathop{\rm sup}\limits_{0\leq s\leq T} |Y^{0,\xi}_{s}-Y^{0,\xi'}_{s} |^{p} + (\int_{0}^{T}
|Z^{0,\xi}_{s}-Z^{0,\xi'}_{s}|^{2}\mathrm{d}s)^{\frac{p}{2}}]\leq {C}_2|\xi-\xi'|^{p}.& \\
\end{array}
\]
\end{corollary}

\noindent \textbf{Proof}: For $t=0$, based on Assumptions $( \mathbf{H2.1}),\ (\mathbf{H2.2})$, and
$$
\left|Y_s^{0,\xi}-Y_s^{0,\xi'}\right|\leq {C_1} \left| X_s^{0,\xi}- X_s^{0,\xi'}\right|,\quad s\in [0,T],
$$
 using Lemma \ref{lem-2}, FBSDE (\ref{equ2.1}) admits a unique $L^p$-solution in the interval $[0,\delta]$. Then, by the inductive method, FBSDE (\ref{equ2.1}) admits a unique $L^p$-solution in the interval $[i\delta,(i+1)\delta],\ 0\leq i\leq k$, where $T=k\delta$ (without loss of generality).

Similar to the proof in Theorem \ref{the-1}, we can extend the $L^P$ estimations from the interval $[i\delta,(i+1)\delta]$ to $[0,T]$. Thus, FBSDE (\ref{equ2.1}) admits a unique $L^p$-solution in interval $[0,T]$ . The proof is complete. $ \qquad \qquad \Box$

\section{Linear quadratic optimal control problem}\label{sec:3}

In the following, we use Theorem \ref{the-1} to study a linear-quadratic optimal control problem. The controlled stochastic system is expressed as follows:
 \begin{equation}
\label{equ3.1}\left \{
\begin{array}
[c]{llll}%
\mathrm{d}X^{u}_s & = &\left [A_sX^{u}_s+B_su_s+b_s\right]\mathrm{d}s + \left[C_sX^u_s+D_su_s+\sigma_s\right] \mathrm{d}B_{s}, & \\
\mathrm{d}X^u_t & = &x,
\end{array}
\right.
\end{equation}
where $A,B,C$, and $D$ are given bounded stochastic matrix-valued functions with proper dimensions, and $b,\sigma$ are vector-valued progressively measurable processes. $X^u(\cdot):[t,T]\times \Omega\to \mathbb{R}^n$, and the set of controls:
$$
\mathcal{U}[t,T]=\{u(\cdot):[t,T]\times \Omega\to \mathbb{R}^m |\ u(\cdot)\ \text{is quadratic integrable progressively measurable process} \}.
$$
The cost functional is given as follows:
\begin{equation*}%
\begin{array}
[c]{ll}%
J(t,x;u(\cdot))=&E\big[\int_t^T\left(\langle Q_sX^u_s,X^u_s \rangle+2\langle S_sX^u_s,u_s \rangle  +\langle R_su_s,u_s \rangle +2\langle q_s, X^u_s \rangle+2\langle \rho_s,u_s\rangle \right)\mathrm{d}s \\
&+ \langle H X^u_T,X^u_T \rangle+2\langle h,X^u_T \rangle\big]. \\
\end{array}
\end{equation*}
Thus, the related optimal control problem is to find an optimal control $\bar{u}(\cdot)\in \mathcal{U}[t,T]$ such that
$$
J(t,x;\bar{u}(\cdot))=\inf_{u\in \mathcal{U}[t,T]}J(t,x;u(\cdot)).
$$

We assume that $(\bar{u}(\cdot),\bar{X}(\cdot))$ is an optimal pair of optimal control problems with a cost functional $J(t,x;u(\cdot))$. For any given control $u(\cdot)\in \mathcal{U}[t,T]$, it follows that
$$
{X}^{\bar{u}+\varepsilon u}_s=\bar{X}_s+\varepsilon X^{0,u}_s,\ t\leq s\leq T,
$$
where ${X}^{\bar{u}+\varepsilon u}(\cdot)$ is the solution of Equation (\ref{equ3.1}) with control $\bar{u}(\cdot)+\varepsilon u(\cdot)$ and $X^{0,u}(\cdot)$ is the solution of the following equation
 \begin{equation}
\left \{
\begin{array}
[c]{llll}%
\mathrm{d}X^{0,u}_s & = &\left [A_sX^{0,u}_s+B_su_s\right]\mathrm{d}s + \left[C_sX^{0,u}_s+D_su_s\right] \mathrm{d}B_{s}, & \\
\mathrm{d}X^{0,u}_t & = &0.
\end{array}
\right.
\end{equation}

Thus, from
$$
0=\lim_{\varepsilon\to 0}\frac{J(t,x;\bar{u}(\cdot)+\varepsilon u(\cdot))-J(t,x;\bar{u}(\cdot))}{\varepsilon},
$$
we have that
$$
0=E\big[\int_t^T \big(\langle Q_s\bar{X}_s+S^{\top}_s\bar{u}_s+q_s,X^{0,u}_s \rangle
+\langle R_s\bar{u}_s+S_s\bar{X}_s+\rho_s,u_s \rangle \big)\mathrm{d}s+\langle H\bar{X}_T+h,X^{0,u}_T \rangle\big].
$$

Then, using the following adjoint equation
\begin{equation}
\label{equ3.2}\left \{
\begin{array}
[c]{llll}%
dY_{s} & = & -\big[A_sY_s+C_sZ_s+Q_s\bar{X}_s+S_s^{\top}\bar{u}_s+q_s  \big]\mathrm{d}s+Z_s\mathrm{d}B_s, & \\
Y_{T} & = & H\bar{X}_T+h, &
\end{array}
\right.
\end{equation}
it follows that for any $u(\cdot)\in \mathcal{U}[t,T]$,
$$
B^{\top}_sY_s+D^{\top}_sZ_s+S_s\bar{X}_s+R_s\bar{u}_s+\rho_s=0,\quad a.e.\ s\in [t,T].
$$
Thus when $R$ is positive, we have
$$
\bar{u}_s=-R^{-1}_s\big[ B^{\top}_sY_s+D^{\top}_sZ_s+S_s\bar{X}_s+\rho_s \big].
$$

We put $\bar{u}_s$ into controlled stochastic system (\ref{equ3.1}) and adjoint Equation (\ref{equ3.2}),
\begin{equation}
\label{equ3.3}\left \{
\begin{array}
[c]{llll}%
\mathrm{d}\bar{X}_s  = &\left [(A_s-B_sR^{-1}_sS_s)\bar{X}_s-B_sR^{-1}_sB^{\top}_sY_s-B_sR^{-1}_s D^{\top}_sZ_s-B_sR^{-1}_s\rho_s+b_s\right]\mathrm{d}s \\
&+\left[(C_s-D_sR^{-1}_sS_s)\bar{X}_s-D_sR^{-1}_sB^{\top}_sY_s-D_sR^{-1}_s D^{\top}_sZ_s-D_sR^{-1}_s\rho_s +\sigma_s\right] \mathrm{d}B_{s}, & \\
dY_{s} = & -\big[(Q_s-S_s^{\top}R^{-1}_sS_s)\bar{X}_s+(A_s-S_s^{\top}R^{-1}_sB^{\top}_s)Y_s
+(C_s-S_s^{\top}R^{-1}_sD^{\top}_s)Z_s-S_s^{\top}R^{-1}_s\rho_s+q_s  \big]\mathrm{d}s\\
&+Z_s\mathrm{d}B_s, & \\
Y_{T}=& H\bar{X}_T+h,\ \bar{X}_t  = x.  &
\end{array}
\right.
\end{equation}

To guarantee the solvability of FBSDE (\ref{equ3.3}), we add the following assumptions.
\begin{description}
\item[$( \mathbf{H3.1})$.] $A_s,B_s,C_s,D_s,\ s\geq 0$ are bounded random matrices.
\end{description}
\begin{description}
\item[$( \mathbf{H3.2})$.] $Q_s-S_s^{\top}R^{-1}_sS_s\geq 0$, $R_s>\delta I,\ \delta>0,\ s\geq 0$, and $H\geq 0$.
\end{description}
\begin{description}
\item[$( \mathbf{H3.3})$.] The norm of matrix $D_s$ is sufficiently small, that is, $|D_s|=\sqrt{\text{tr}(D_sD_s^{\top})},\ s\geq 0$.
\end{description}

\begin{lemma}
\label{lem-3-1} Let Assumptions $( \mathbf{H3.1}),\ (\mathbf{H3.2})$
hold, FBSDE (\ref{equ3.3}) admits a unique $L^2$-solution $(\bar{X},Y,Z)$. For the initial value of $\bar{X}_t=\xi$ and $\bar{X}'_t=\xi'$, where  $\xi,\xi' \in L^{p}(\Omega,\mathcal{F}_{t}, P;\mathbb{R}^{n})$. Then, we can obtain the $L^2$ estimations of FBSDE (\ref{equ3.3}),
\begin{equation}
\label{equa3.3-1}
E [ \mathop{\rm sup}\limits_{t\leq s\leq T}
|\bar{X}_{s}-\bar{X}'_{s} |^{2} + \mathop{\rm sup}\limits_{t\leq s\leq T} |Y_{s}-Y'_{s} |^{2} + \int_{t}^{T}
|Z_{s}-Z'_{s}|^{2}\mathrm{d}s\mid \mathcal{F}_{t} ]\leq {C}|\xi-\xi'|^{2},
\end{equation}
where $(\bar{X},Y,Z)$ is the solution of (\ref{equ3.3}) with the initial value $\xi$,  $(\bar{X}',Y',Z')$ with $\xi'$, and $C$ depends on the coefficients of FBSDE (\ref{equ3.3}).
\end{lemma}
\noindent \textbf{Proof}:
Denoting the coefficients of the fully coupled FBSDE (\ref{equ3.3}) as follows:
\begin{equation}\label{equ3.3-0}
F(s,x,y,z)=
\begin{pmatrix}
-(Q_s-S_s^{\top}R^{-1}_sS_s)x-(A_s-S_s^{\top}R^{-1}_sB^{\top}_s)y
-(C_s-S_s^{\top}R^{-1}_sD^{\top}_s)z\\
(A_s-B_sR^{-1}_sS_s)x-B_sR^{-1}_sB^{\top}_sy-B_sR^{-1}_s D^{\top}_sz\\
(C_s-D_sR^{-1}_sS_s)x-D_sR^{-1}_sB^{\top}_sy-D_sR^{-1}_s D^{\top}_sz\\
\end{pmatrix}
.
\end{equation}
Based on Assumption $( \mathbf{H3.2})$, one obtains
$$
\langle F(s,x,y,z),(x,y,z)  \rangle=-\left[ \langle (Q_s-S_s^{\top}R^{-1}_sS_s)x,x \rangle
+\langle R_s^{-1}(B_s^{\top}y+D_s^{\top}z),(B_s^{\top}y+D_s^{\top}z) \rangle  \right].
$$

Thus, there exist constants $c_1\geq 0$ and $c_2>0$ such that
$$
\langle F(s,x,y,z),(x,y,z)  \rangle\leq -c_1|x|^2-c_2| B_s^{\top}y+D_s^{\top}z |^2.
$$
and
$$
\langle H x,x \rangle\geq 0.
$$
These are the monotonicity conditions in \citet{HP95} and \citet{PW99}. Then, by Theorem 3.1 of \citet{PW99}, FBSDE (\ref{equ3.3}) admits a unique $L^2$-solution $(\bar{X},Y,Z)$.

Now, we calculate the distance between $(\bar{X},Y,Z)$ and $(\bar{X}',Y',Z')$ with initial values $(\xi,\xi')$. First, by Assumption $( \mathbf{H3.1})$, for the forward SDE of $\bar{X}$ and $\bar{X}'$, we have the inequality of $\hat{X}=\bar{X}-\bar{X}'$,
\begin{eqnarray}
\begin{split}\label{eq:3.3-2}
& E_t\bigg[\displaystyle\sup_{t\leqslant s\leqslant
T}|\hat{X}_s|^2\bigg]+
 E_t\bigg[\displaystyle\int_{t}^T|\hat X_s|^2\mathrm{d}s\bigg]\leqslant K\bigg\{ E_t\bigg[|\xi-\xi'|^2\bigg]+
 E_t\bigg[\int_{t}^T(|B_s^{\top}\hat{Y}_s+D_s^{\top}\hat Z_s|^2)\mathrm{d}s\bigg]\bigg\},
\end{split}
\end{eqnarray}
where $E_t[\cdot]$ is the conditional expectation based on information $\mathcal{F}_t$. Second, from Assumption $( \mathbf{H3.1})$, for the backward SDE of $(Y,Z)$ and $(Y',Z')$, it follows that
\begin{eqnarray}
    \begin{split}\label{eq:3.3-3}
& E_t\bigg[\displaystyle\sup_{t\leqslant s\leqslant
T}|\hat{Y}_s|^2\bigg]+E_t\bigg[\displaystyle\int_{t}^T(|\hat Y_s|^2+|\hat Z_s|^2)\mathrm{d}s\bigg]\leqslant K\bigg\{ E_t\bigg[\int_t^T|\hat{X}_s|^2\mathrm{d}s\bigg]+
 E_t\big[| \hat{X}_T|^2\big]\bigg\}.
\end{split}
  \end{eqnarray}

Applying It\^{o} formula to $\langle\hat X (s),\hat Y (s)\rangle,$ we have
\begin{eqnarray}\label{eq:3.3-4}
  \begin{split}
 &E_t[\langle\hat{X}_T, H\hat{X}_T\rangle]
=E_t[\int_t^T \langle F(s, U_s)-F(s, U'_s),
U_s-U'_s\rangle\mathrm{d}s]+E_t[\langle\hat{Y}_t, \xi-\xi'\rangle],
\end{split}
\end{eqnarray}
where $U_s=(\bar{X}_s,Y_s,Z_s)$, $U'_s=(\bar{X}'_s,Y'_s,Z'_s)$, and $F(\cdot)$ is given in (\ref{equ3.3-0}). From the monotone properties of $F(\cdot)$ and $H$, we have
\begin{eqnarray}\label{eq:3.3-5}
  \begin{split}
    &
 E_t\bigg[\int_{t}^T(|B_s^{\top}\hat{Y}_s+D_s^{\top}\hat Z_s|^2)\mathrm{d}s\bigg]\leq
 \frac{\varepsilon}{c_2} E_t\bigg[\displaystyle\sup_{t\leqslant s\leqslant
T}|\hat{Y}_s|^2\bigg]+\frac{1}{4\varepsilon c_2}E_t[|\xi-\xi'|^2],
  \end{split}
\end{eqnarray}
where $\varepsilon>0$, which is given later. Combining inequalities (\ref{eq:3.3-2}) and (\ref{eq:3.3-3}), it follows that
\begin{eqnarray}
    \begin{split}\label{eq:3.3-6}
& E_t\bigg[\displaystyle\sup_{t\leqslant s\leqslant
T}|\hat{Y}_s|^2\bigg]+E_t\bigg[\displaystyle\int_{t}^T(|\hat Y_s|^2+|\hat Z_s|^2)\mathrm{d}s\bigg]\\
\leqslant &K^2\bigg\{ E_t\bigg[|\xi-\xi'|^2\bigg]+
 E_t\bigg[\int_{t}^T(|B_s^{\top}\hat{Y}_s+D_s^{\top}\hat Z_s|^2)\mathrm{d}s\bigg]\bigg\}.
\end{split}
  \end{eqnarray}
Now, applying inequality (\ref{eq:3.3-5}), we can obtain that
\begin{eqnarray}
    \begin{split}\label{eq:3.3-7}
& E_t\bigg[\displaystyle\sup_{t\leqslant s\leqslant
T}|\hat{Y}_s|^2\bigg]+E_t\bigg[\displaystyle\int_{t}^T(|\hat Y_s|^2+|\hat Z_s|^2)\mathrm{d}s\bigg]\\
\leqslant &
\frac{\varepsilon K^2}{c_2} E_t\bigg[\displaystyle\sup_{t\leqslant s\leqslant
T}|\hat{Y}_s|^2\bigg]+
\frac{(1+4\varepsilon c_2)K^2}{4\varepsilon c_2}E_t[|\xi-\xi'|^2].
\end{split}
  \end{eqnarray}
Let $\displaystyle \varepsilon<\frac{c_2 }{K^2}$, then there exists a constant $C>0$ such that
\begin{eqnarray}
    \begin{split}\label{eq:3.3-8}
& E_t\bigg[\displaystyle\sup_{t\leqslant s\leqslant
T}|\hat{Y}_s|^2\bigg]+E_t\bigg[\displaystyle\int_{t}^T(|\hat Y_s|^2+|\hat Z_s|^2)\mathrm{d}s\bigg]
\leqslant C E_t[|\xi-\xi'|^2],
\end{split}
  \end{eqnarray}
where $C$ depends on the coefficients of FBSDE (\ref{equ3.3}). Then, by combining (\ref{eq:3.3-2}) and (\ref{eq:3.3-5}), we can obtain inequality (\ref{equa3.3-1}). The proof is complete. $ \qquad \qquad \Box$

\bigskip

The main results of this section are given as follows:
\begin{theorem}
\label{the-4} Let Assumptions $( \mathbf{H3.1}),\ (\mathbf{H3.2}), \ (\mathbf{H3.3})$
hold. Then, FBSDE (\ref{equ3.3}) admits a unique $L^p$-solution $(\bar{X},Y,Z)$ with $p>2$.
\end{theorem}
\noindent \textbf{Proof}: Based on Assumptions $( \mathbf{H3.1}),\ (\mathbf{H3.2})$ and Lemma \ref{lem-3-1}, we have that FBSDE (\ref{equ3.3}) admits a unique $L^2$-solution $(\bar{X},Y,Z)$. Let $(\bar{X},Y,Z)$ be the solution of (\ref{equ3.3}) with initial value $\xi$,  $(\bar{X}',Y',Z')$ with $\xi'$, where $\xi,\xi' \in L^{p}(\Omega,\mathcal{F}_{t}, P;\mathbb{R}^{n})$. By using Lemma \ref{lem-3-1}, we have
\begin{equation}
\label{equa3.4-1}
E [ \mathop{\rm sup}\limits_{t\leq s\leq T}
|\bar{X}_{s}-\bar{X}'_{s} |^{2} + \mathop{\rm sup}\limits_{t\leq s\leq T} |Y_{s}-Y'_{s} |^{2} + \int_{t}^{T}
|Z_{s}-Z'_{s}|^{2}\mathrm{d}s\mid \mathcal{F}_{t} ]\leq {C}|\xi-\xi'|^{2}.
\end{equation}
Then, combining Assumption $( \mathbf{H3.3})$ and Theorem \ref{the-1}, FBSDE (\ref{equ3.3}) admits a unique $L^p$-solution $(\bar{X},Y,Z)$. The proof is complete. $ \qquad \qquad \Box$

\bigskip

In the following, we extend the results of Theorem \ref{the-4} to a general linear coupled FBSDE,
\begin{equation}
\label{equ3.5}\left \{
\begin{array}
[c]{llll}%
\mathrm{d}X_{s} & = & [A_1X_s+B_1Y_s+D_1Z_s]\mathrm{d}s + [A_2X_s+B_2Y_s+D_2Z_s] \mathrm{d}B_{s}, & \\
\mathrm{d}Y_{s} & = & -[A_3X_s+B_3Y_s+D_3Z_s]\mathrm{d}s +
Z_{s}\mathrm{d}B_{s},  & \\
Y_{T} & = & HX_T, \ X(0)=x,&
\end{array}
\right.
\end{equation}
where $(A_i,B_i,D_i)_{i=1}^3$ are bounded stochastic matrices of time. For convenience, we omit time $s$.
\begin{corollary}
\label{coro-4} Let the coefficients of FBSDE (\ref{equ3.5}) satisfy the monotonicity conditions, and the norm of $D_2$ is sufficiently small. Then, FBSDE (\ref{equ3.5}) admits a unique $L^p$-solution $({X},Y,Z)$ with $p>2$.
\end{corollary}
\noindent \textbf{Proof}: The monotonicity conditions of FBSDE (\ref{equ3.5}) show that  FBSDE (\ref{equ3.3}) admits a unique $L^2$-solution with $L^2$ estimations, which is the same as the results given in Lemma \ref{lem-3-1}. Thus, FBSDE (\ref{equ3.5}) admits a unique $L^p$-solution $({X},Y,Z)$. The proof is complete. $ \qquad \qquad \Box$

\section{Special case of FBSDEs }\label{sec:4}

Now, we consider the following FBSDE:
\begin{equation}
\label{equ4.1}\left \{
\begin{array}
[c]{llll}%
dX_{s}^{0,\xi} & = & b(s,X_{s}^{0,\xi},Y_{s}^{0,\xi},Z_{s}^{0,\xi})\mathrm{d}s + \sigma
(s,X_{s}^{0,\xi},Y_{s}^{0,\xi}) \mathrm{d}B_{s}, & \\
dY_{s}^{0,\xi} & = & -f(s,X_{s}^{0,\xi},Y_{s}^{0,\xi},Z_{s}^{0,\xi})\mathrm{d}s +
Z_{s}^{0,\xi}\mathrm{d}B_{s}, \  \  \  \  \ s\in[0,T], & \\
X_{0}^{0,\xi} & = & \xi, & \\
Y_{T}^{0,\xi} & = & \Phi(X_{T}^{0,\xi}), &
\end{array}
\right.
\end{equation}
where the diffusion term $\sigma(s,\cdot)$ in forward SDE does not depend on $Z_{s}^{0,\xi},\ s\geq 0$.

We first introduce the results of \citet{CZ13}.
\begin{lemma}\label{lem-4}
Let Assumptions $( \mathbf{H2.1})$ and $( \mathbf{H2.2})$ hold, there exists a random field $u(t,x)$ satisfying

(i). $u(T,x)=\Phi(x)$;

(ii). $u(t,x)$ is $\mathcal{F}_t$ measurable;

(iii). $|u(t,x_1)-u(t,x_2)|\leq K|x_1-x_2 |,\ x_1,x_2\in \mathbb{R}^n$;

(iv). For any $0\leq t_1\leq t_2\leq T$ such that $|t_2-t_1|\leq \delta(K)$, where $\delta(K)$ is a sufficiently small constant, a unique solution to FBSDE (\ref{equ4.1}) over $[t_1,t_2]$, satisfies $Y^{0,\xi}_{t_2}=u(t_2,X^{0,\xi}_{t_2})$ and
 $Y^{0,\xi}_{t_1}=u(t_1,X^{0,\xi}_{t_1})$.
\bigskip

Thus, for any given $T>0$, FBSDE (\ref{equ4.1}) admits a unique solution in the interval $[0,T]$ and
$Y^{0,\xi}_{t}=u(t,X_t^{0,\xi})$.
\end{lemma}

Based on Lemma \ref{lem-4}, we can obtain $L^2$ estimations for the solution of FBSDE (\ref{equ4.1}).
\begin{lemma}\label{lem-4-1}
Let Assumptions $( \mathbf{H2.1})$ and $( \mathbf{H2.2})$ hold, and $u(t,x)$ satisfies the conditions in Lemma \ref{lem-4}. Then, we obtain the $L^2$ estimations of FBSDE (\ref{equ4.1}),
\begin{equation}
\label{equ4.2}
E [ \mathop{\rm sup}\limits_{0\leq s\leq T}
|X^{0,\xi}_{s}-X^{0,\xi'}_{s} |^{2} + \mathop{\rm sup}\limits_{0\leq s\leq T} |Y^{0,\xi}_{s}-Y^{0,\xi'}_{s} |^{2} + \int_{0}^{T}
|Z^{0,\xi}_{s}-Z^{0,\xi'}_{s}|^{2}\mathrm{d}s\mid \mathcal{F}_{0} ]\leq {C}|\xi-\xi'|^{2},
\end{equation}
and thus
\begin{equation}
\label{equ4.3}
E [ \mathop{\rm sup}\limits_{t\leq s\leq T}
|X^{0,\xi}_{s}-X^{0,\xi'}_{s} |^{2} + \mathop{\rm sup}\limits_{t\leq s\leq T} |Y^{0,\xi}_{s}-Y^{0,\xi'}_{s} |^{2} + \int_{t}^{T}
|Z^{0,\xi}_{s}-Z^{0,\xi'}_{s}|^{2}\mathrm{d}s\mid \mathcal{F}_{t} ]\leq {C}|X^{0,\xi}_{t}-X^{0,\xi'}_{t} |^{2} ,
\end{equation}
where $C>0$ is a constant which depends on the constants $L$ and $K$ given in $( \mathbf{H2.1})$ and $( \mathbf{H2.2})$.
\end{lemma}

\noindent \textbf{Proof}: For a given partition $0=t_0<t_1< \cdots<t_k=T$, satisfies $|t_i-t_{i-1}|<\delta(K)$. Thus, we have the $L^2$ estimations for $1\leq i\leq k$,
$$
E [ \mathop{\rm sup}\limits_{t_{i-1}\leq s\leq t_i}
|X^{0,\xi}_{s}-X^{0,\xi'}_{s} |^{2} + \mathop{\rm sup}\limits_{t_{i-1}\leq s\leq t_i} |Y^{0,\xi}_{s}-Y^{0,\xi'}_{s} |^{2} + \int_{t_{i-1}}^{t_i}
|Z^{0,\xi}_{s}-Z^{0,\xi'}_{s}|^{2}\mathrm{d}s\mid \mathcal{F}_{t_{i-1}} ]\leq {C}_0|X^{0,\xi}_{t_{i-1}}-X^{0,\xi'}_{t_{i-1}}|^{2}.
$$

First, we consider the case $i=1$,
$$
E [ \mathop{\rm sup}\limits_{t_{0}\leq s\leq t_1}
|X^{0,\xi}_{s}-X^{0,\xi'}_{s} |^{2} + \mathop{\rm sup}\limits_{t_{0}\leq s\leq t_1} |Y^{0,\xi}_{s}-Y^{0,\xi'}_{s} |^{2} + \int_{t_{0}}^{t_1}
|Z^{0,\xi}_{s}-Z^{0,\xi'}_{s}|^{2}\mathrm{d}s\mid \mathcal{F}_{t_{0}} ]\leq {C}_0|\xi-\xi'|^{2}
$$
and $i=2$,
$$
E [ \mathop{\rm sup}\limits_{t_{1}\leq s\leq t_2}
|X^{0,\xi}_{s}-X^{0,\xi'}_{s} |^{2} + \mathop{\rm sup}\limits_{t_{1}\leq s\leq t_2} |Y^{0,\xi}_{s}-Y^{0,\xi'}_{s} |^{2} + \int_{t_{1}}^{t_2}
|Z^{0,\xi}_{s}-Z^{0,\xi'}_{s}|^{2}\mathrm{d}s\mid \mathcal{F}_{t_{1}} ]\leq {C}_0|X^{0,\xi}_{t_{1}}-X^{0,\xi'}_{t_{1}}|^{2}.
$$

Based on a similar idea in the proof of Theorem \ref{the-1}, there exists a constant $C_1>0$ such that
$$
E [ \mathop{\rm sup}\limits_{t_{0}\leq s\leq t_2}
|X^{0,\xi}_{s}-X^{0,\xi'}_{s} |^{2} + \mathop{\rm sup}\limits_{t_{0}\leq s\leq t_2} |Y^{0,\xi}_{s}-Y^{0,\xi'}_{s} |^{2} + \int_{t_{0}}^{t_2}
|Z^{0,\xi}_{s}-Z^{0,\xi'}_{s}|^{2}\mathrm{d}s\mid \mathcal{F}_{t_{0}} ]\leq {C}_1|\xi-\xi'|^{2}.
$$
Then, using the inductive method,  we can obtain inequality (\ref{equ4.2}). In a similar manner to the proof in inequality (\ref{equ4.2}), we can establish inequality (\ref{equ4.3}), which completes this proof. $ \qquad \qquad \Box$

\begin{theorem}\label{the-5}
Let the conditions in Lemma \ref{lem-4} hold. Then, FBSDE (\ref{equ4.1}) admits a unique $L^p$-solution with $p>2$.
\end{theorem}
\noindent \textbf{Proof}: Combining Lemma \ref{lem-4} and Lemma \ref{lem-4-1}, FBSDE (\ref{equ4.1}) admits a unique $L^2$-solution. Applying Theorem \ref{the-1}, FBSDE (\ref{equ4.1}) admits a unique $L^p$-solution. The proof is complete. $ \qquad \qquad \Box$

\bigskip
Based on a method similar to the proof of Theorem \ref{the-5}, we can improve the $L^2$-solution of Theorem 11.3.3 in \citet{CZ13} and  Theorem 7.3 in \citet{MW15} to the $L^p$-solution with $p>2$.
\begin{corollary}\label{coro-5}
Let Assumptions $( \mathbf{H2.1})$ and $( \mathbf{H2.2})$ hold, $b,\sigma,f,\Phi$ are deterministic functions, and $\sigma>\delta I$. Then, FBSDE (\ref{equ4.1}) admits a unique $L^p$-solution with $p>2$.
\end{corollary}

\begin{corollary}\label{coro-6}
Let the conditions of Theorem 7.3 in \citet{MW15} be correct. Then, FBSDE (\ref{equ4.1}) admits a unique $L^p$-solution with $p>2$.
\end{corollary}

\section{Extensions of the main results}\label{sec:4.5}
In this section, we apply the results of \citet{Y20} to improve the main results of Theorem \ref{the-1}. We first introduce a constant $K_p,\ p>1$ which is given by \citet{Y20},
$$
K_p=\overline{K}^{1/p}_p\big(\frac{p}{p+1}+2\underline{K}^{-1/p}_p\frac{2p-1}{p-1}\big),
$$
where $\underline{K}_p,\overline{K}_p$ satisfy the following Burkholder-Davis-Gundy's inequalities,
$$
\underline{K}_pE_t\bigg(\int_t^T\left|Z_s\right|^2\mathrm{d}s\bigg)^{p/2}\leq E_t\bigg(\sup_{r\in[t,T]}\left|\int_t^rZ_s\mathrm{d}B_s\right|^p\bigg)\leq \overline{K}_pE_t\bigg(\int_t^T\left|Z_s\right|^2\mathrm{d}s\bigg)^{p/2}.
$$

Based on constant $K_p$, we give the following assumption.
\begin{description}
\item[$( \mathbf{H5.1})$] There exist constants $K$ and $L_{\sigma}\geq0$ such that for all $t\in[0,T],\ x_{1},x_{2}%
\in \mathbb{R}^{n},\ y_{1},y_{2}\in \mathbb{R}^m,\ z_{1},z_{2}\in \mathbb{R}%
^{m}$,
\[
|b(t,x_{1},y_{1},z_{1})-b(t,x_{2},y_{2},z_{2})|\leq
K(|x_{1}-x_{2}|+|y_{1}-y_{2}|+|z_{1}-z_{2}|),
\]
\[
|\sigma(t,x_{1},y_{1},z_{1})-\sigma(t,x_{2},y_{2},z_{2})|\leq
K(|x_{1}-x_{2}|+|y_{1}-y_{2}|)+L_{\sigma}|z_{1}-z_{2}|,
\]
\[
|f(t,x_{1},y_{1},z_{1})-f(t,x_{2},y_{2},z_{2})|\leq
K(|x_{1}-x_{2}|+|y_{1}-y_{2}|+|z_{1}-z_{2}|),
\]
\[
|\Phi(x_{1})-\Phi(x_{2})|\leq
K|x_{1}-x_{2}|,
\]
\[
K_pL_{\sigma}K<1.
\]
\end{description}

We introduce the results of Theorem 2.3 of \citet{Y20} as follows. More details see \citet{Y18}, in which they further established the probabilistic interpretation for a system of quasilinear parabolic partial differential-algebraic equations by fully coupled FBSDEs.
\begin{lemma}
\label{lem-5-1} Let Assumptions $( \mathbf{H2.1}),\ (\mathbf{H5.1})$
hold. Then, for any given $p> 2,$ there exist constants ${\delta} > 0$ depending on $(K,L_\sigma)$, and ${C}_0$ depends on $(p,L,K,L_\sigma)$, such that
for every $\xi,\xi' \in L^{p}(\Omega,\mathcal{F}_{t}, P;\mathbb{R}^{n}),$
\[%
\begin{array}
[c]{llll}%
 &  E [ \mathop{\rm sup}\limits_{t\leq s\leq t+\delta}
|X^{t,\xi}_{s} |^{p} + \mathop{\rm sup}\limits_{t\leq s\leq
t+\delta} |Y^{t,\xi}_{s} |^{p} + (\int_{t}^{t+\delta}
|Z^{t,\xi}_{s}|^{2}\mathrm{d}s)^{\frac{p}{2}} \mid \mathcal{F}_{t}]\leq {C}_0(1 + |\xi|^{p}).& \\
\end{array}
\]
Then, FBSDE (\ref{equ2.1}) admits a unique $L^p$ solution in the interval $[t,t+\delta]$ with $p>2$.
\end{lemma}

Based on Lemma \ref{lem-5-1}, we show the relation between $L^2$ and $L^p$ estimations under Assumptions $( \mathbf{H2.1}),\ (\mathbf{H5.1})$.
\begin{theorem}
\label{the-5-1} Let Assumptions $( \mathbf{H2.1}),\ (\mathbf{H5.1})$
hold, and we assume that for every $\xi,\xi' \in L^{p}(\Omega,\mathcal{F}_{t}, P;\mathbb{R}^{n})$, $L^2$ estimations of FBSDE (\ref{equ2.1})  are right,
\begin{equation}
 E [ \mathop{\rm sup}\limits_{t\leq s\leq T}
|X^{t,\xi}_{s} |^{2} + \mathop{\rm sup}\limits_{t\leq s\leq T} |Y^{t,\xi}_{s} |^{2} + (\int_{t}^{T}
|Z^{t,\xi}_{s}|^{2}\mathrm{d}s)\mid \mathcal{F}_{t} ]\leq {C}_1(1 + |\xi|^{2}),
\end{equation}
\begin{equation}
E [ \mathop{\rm sup}\limits_{t\leq s\leq T}
|X^{t,\xi}_{s}-X^{t,\xi'}_{s} |^{2} + \mathop{\rm sup}\limits_{t\leq s\leq T} |Y^{t,\xi}_{s}-Y^{t,\xi'}_{s} |^{2} + \int_{t}^{T}
|Z^{t,\xi}_{s}-Z^{t,\xi'}_{s}|^{2}\mathrm{d}s\mid \mathcal{F}_{t} ]\leq {C}_1|\xi-\xi'|^{2},
\end{equation}
where $C_1$ is a constants and independent from $t\in [0,T]$. Furthermore, we assume that $K_pL_{\sigma}\sqrt{C_1}<1$.

Then, FBSDE (\ref{equ2.1}) admits a unique $L^p$ ($p>2$)
solution with $t=0$, and
\[%
\begin{array}
[c]{llll}%
 &  E [ \mathop{\rm sup}\limits_{0\leq s\leq T}
|X^{0,\xi}_{s} |^{p} + \mathop{\rm sup}\limits_{0\leq s\leq T} |Y^{0,\xi}_{s} |^{p} + (\int_{0}^{T}
|Z^{0,\xi}_{s}|^{2}\mathrm{d}s)^{\frac{p}{2}}]\leq {C}_2(1 + |\xi|^{p}),& \\
 &  E [ \mathop{\rm sup}\limits_{0\leq s\leq T}
|X^{0,\xi}_{s}-X^{0,\xi'}_{s} |^{p} + \mathop{\rm sup}\limits_{0\leq s\leq T} |Y^{0,\xi}_{s}-Y^{0,\xi'}_{s} |^{p} + (\int_{0}^{T}
|Z^{0,\xi}_{s}-Z^{0,\xi'}_{s}|^{2}\mathrm{d}s)^{\frac{p}{2}}]\leq {C}_2|\xi-\xi'|^{p}.& \\
\end{array}
\]
\end{theorem}
\noindent \textbf{Proof}: Based on Assumptions $( \mathbf{H2.1}),\ (\mathbf{H5.1})$ and $K_pL_{\sigma}\sqrt{C_1}<1$, using Lemma \ref{lem-5-1}, we can show that there exists $\delta>0$ such that
\[%
\begin{array}
[c]{llll}%
 & & E[ \mathop{\rm sup}\limits_{(i-1)\delta\leq s\leq i\delta}
|X^{0,\xi}_{s} |^{p}+ \mathop{\rm sup}\limits_{(i-1)\delta\leq s\leq i\delta}|Y^{0,\xi}_{s} |^{p}+ (\int_{(i-1)\delta}^{i\delta}
|Z^{0,\xi}_{s}|^{2}ds)^{\frac{p}{2}} \mid \mathcal{F}_{(i-1)\delta}]\leq \hat{C}_1(1 + |X^{0,\xi}_{(i-1)\delta}|^{p}), \\
\end{array}
\]
where $1\leq i\leq k$, $T=k\delta$, and $\hat{C}_1$ depends on constants $C_1$ and $L,K,L_{\sigma}$ in Assumptions $(\mathbf{H2.1}),\ (\mathbf{H5.1})$. The following proof is the same as that in Theorem \ref{the-1}. Thus, we complete the proof. $ \qquad \qquad \Box$

\section{Conclusion}\label{sec:5}

We studied whether an adapted $L^2$-solution of fully  coupled FBSDEs is an adapted $L^P$-solution for some $p>2$ which was proposed in \citet{Y20}. Based on the usual Lipschatiz (Lipschatiz constant of $Z$ in diffusion term should be sufficiently small) and linear growth conditions on the coefficients, we established  a uniform $L^p$ estimations for fully coupled FBSDEs in different small time intervals. Then, we extend the small time interval $L^p$ estimations to a global time interval. That is, for a given terminal time $T>0$,  we proved that the unique $L^2$-solution of fully coupled FBSDE is an $L^p$-solution with $p>2$.

Based on the main results of this study, we further considered the fully coupled linear FBSDEs which are generalized by a linear quadratic optimal control problem, and established the $L^p$ estimations for the fully coupled linear FBSDEs with random coefficients.   We also improved the $L^2$-solution of fully coupled FBSDEs to the $L^p$-solution based on the "decoupling random field" method.

\bibliography{gexp1}

\end{document}